\documentclass[11pt,a4paper]{article}
\usepackage[latin2]{inputenc}
\usepackage[english]{babel}
\usepackage{a4wide}
\usepackage[english]{babel}
\usepackage[nottoc]{tocbibind}
\usepackage{amsmath,amsfonts,amstext,amsbsy,amsopn,amssymb,amsthm,amscd}
\usepackage{amsxtra,latexsym}
\usepackage{bbm}
\usepackage{exscale}
\usepackage{paralist}
\usepackage[article]{math_env} %von Thomas Merkle
\numberwithin{equation}{section}
\usepackage{dmath_abbr_wias}
\DeclareMathOperator{\curl}{curl}
\DeclareMathOperator{\sym}{sym}
\usepackage{graphicx,mathrsfs}
\usepackage{psfrag}
\DeclareMathOperator{\skyw}{skew}
\DeclareMathOperator{\Curl}{Curl}
\DeclareMathOperator{\dyw}{div}
\DeclareMathOperator{\Dyw}{Div}
\DeclareMathOperator{\Sym}{Sym}
\DeclareMathOperator{\so}{\mathfrak{so}}

\newcommand{\nn}{\nonumber}

\renewcommand{\D}{{\mathbb C_{{\rm e}}}}
\newcommand{\Dc}{{\mathbb C_{{\rm c}}}}

\newcommand{\Ha}{{\mathbb C_{{\rm micro}}}}
\newcommand{\Lc}{{\mathbb L}_{{\rm c}}}

\newcommand{\di}{{\mathrm d}}

\usepackage{color}
\definecolor{dred}{rgb}{.8,0,0}
\definecolor{ddmagenta}{rgb}{0.7,0,0.9}
\definecolor{ddcyan}{rgb}{0,0.2,1.0}
\definecolor{Green}{rgb}{0.,0.5,0.}

\makeatletter
\let\@fnsymbol\@arabic
\makeatother

\begin{document}
\title{A global higher regularity result for the static relaxed micromorphic model on smooth domains}
\author{Dorothee Knees\thanks{
Institute of Mathematics, University of Kassel, Heinrich-Plett Str. 40, 34132 Kassel, Germany, dknees@mathematik.uni-kassel.de
}\,,
\and
Sebastian Owczarek\thanks{Faculty of Mathematics and Information Science, Warsaw University of Technology, ul. Koszykowa 75, 00-662 Warsaw, Poland, sebastian.owczarek@pw.edu.pl}\,,
\and
Patrizio Neff\thanks{Lehrstuhl f\"ur Nichtlineare Analysis und Modellierung, Fakult\"at f\"ur Mathematik, Universit\"at Duisburg-Essen, Campus Essen, Thea-Leymann Str. 9, 45127 Essen, Germany, patrizio.neff@uni-due.de}
}

\date{ }
%\date{Dorothee: 17. Maerz 2013, 6. Juli 2022}

\maketitle
%% article and other classes: abstract after maketitle

\begin{abstract}
We derive a global higher regularity result for weak solutions of the linear relaxed micromorphic model on smooth domains. The governing equations  consist of a linear elliptic system of partial differential equations that is coupled with a system of Maxwell-type. The result is obtained by combining a Helmholtz decomposition argument  with regularity results for linear elliptic systems and the classical  embedding  of $H(\dyw;\Omega)\cap H_0(\curl;\Omega)$ into $H^1(\Omega)$. 
\end{abstract}

%
%----------------------------------------------------------------------


\noindent
\textbf{AMS Subject Classification 2020:}
%\subjclass[2020]{%
35Q74, %PDEs in connection with mechanics of deformable solids 
35B65, %Smoothness and regularity of solutions to PDEs
49N60, %Regularity of solutions in optimal control
74A35, %Polar materials
74G40.  %Regularity of solutions of equilibrium problems in solid mechanics    
%}
%----------------------------------------------------------------------

%--------------------------------------------------------------------------

%
%------------------------------------------
\section{Introduction}
\label{s:introduction}

The relaxed micromorphic model is a novel generalised continuum model allowing to describe size effects and band-gap behaviour of microstructured solids with effective equations ignoring the detailed microstructure \cite{exprolingalberdi,Neff_unfol_22,Ghibadyn,MADEO2016,NGperspective,Owczghibaneffexist}. As a micromorphic model it couples the classical displacement $u:\Omega\subset\R^3\to\R^3$ with a non-symmetric tensor field $P:\Omega\subset\R^3\to\R^{3\times 3}$ called the microdistortion through the variational problem \begin{equation}
\label{Variational}
\begin{split}
\int_{\Omega}&\Big(\big\langle \D\sym(\nabla u-P),\sym(\nabla u-P)\big\rangle
+\big\langle\Ha\sym P,\sym P\big\rangle
+ \big\langle\Lc\Curl P,\Curl P\big\rangle\\[1ex]
&\hspace{1ex}-\langle f,u\rangle-\langle M,P\rangle\Big)\,\di x\quad\longrightarrow \quad {\rm min\,\,\,\, w.r.t.\,\,\, }(u,P)\,,
\end{split}
\end{equation}
subject to suitable boundary conditions. The tensor $\Lc$ introduces a size-dependence into the model in the sense that smaller samples respond relatively stiffer. The existence and uniqueness in the static case follows from the incompatible Korn's inequality \cite{Lewintan_Korn2020,lewintan_neff_2021,NeffLP,neff2012poincare}. The constitutive tensors $\D$, $\Ha$ and $\Lc$ are to be found by novel homogenisation stra\-tegies \cite{Neffhomog,homo_1,schroder2021lagrange}. Letting $\Ha\rightarrow +\infty$ the models response tends to the linear Cosserat model \cite{Ghiba2022CosseratME}. A range of engineering relevant analytical solutions are already available for the relaxed micromorphic model \cite{rizzi2021torsion}. The solution is naturally found as $u\in H^1(\Omega)$ and $P\in H(\Curl;\Omega)$, thus the microdistortion $P$ may have jumps in normal direction. The implementation in the finite element context needs standard element formulations for the displacement $u$, but e.g. N\'ed\'elec - spaces for $P$ in order to achieve optimal convergence rates \cite{Rizzi_boundary_2022,schroder2021lagrange,SKY2022_primal,sky2021hybrid}.

However, it is sometimes preferred to circumvent the N\'ed\'elec - framework and to work with $H^1(\Omega)$ for the microdistortion tensor $P$. For these cases it is mandatory to clarify in advance whether the regularity of $P$ allows for a faithful result. In this spirit, we continue here the investigation of regularity in the static case and  we will be able to derive a global higher regularity result for weak solutions  of the relaxed micromorphic model (global as opposed to only interior regularity). It extends the local result from \cite{KON23} to smooth domains. The latter is formulated on a bounded domain $\Omega\subset \R^3$ and the Euler-Lagrange equations to \eqref{Variational} read as follows (\cite{exprolingalberdi,Ghibadyn,NGperspective,Owczghibaneffexist,OGNdynamicreg}): given positive definite and symmetric material dependent coefficient tensors $\bbC_e:\Omega\rightarrow \Lin(\Sym(3),\Sym(3))$, $\bbC_\text{micro}:\Omega\rightarrow \Lin(\Sym(3),\Sym(3))$ and $\bbL_c:\Omega\rightarrow \Lin(\R^{3\times 3},\R^{3\times 3})$ determine a displacement field $u:\Omega\to\R^3$ and a non-symmetric microdistortion tensor $P:\Omega\to\R^{3\times 3}$ satisfying 
\begin{align}
\label{Main}
 0&=\Dyw\Big(\bbC_e\sym(\nabla u - P)\Big) + f\quad \text{in }\Omega,\nn\\
 0&= -\Curl\Big(\bbL_c\Curl P\big) + \bbC_e\sym(\nabla u - P) - \bbC_\text{micro}\sym P + M\quad \text{in }\Omega\,
\end{align}
together with suitable boundary conditions. 
Here, $f:\Omega\to\R^3$ is a given volume force density, $M:\Omega\to \R^{3\times 3}$ a body moment tensor and $\sigma=\bbC_e\sym(\nabla u - P)$ is the symmetric force stress tensor while $m=\Lc\Curl P$ is the non-symmetric moment tensor.

The main result of the present contribution is Theorem \ref{thm:globreg} which states that on smooth domains and with smooth coefficients weak solutions of \eqref{Main} are more regular and satisfy 
\begin{equation}
\label{mainreg}
 u\in H^2(\Omega)\,,\quad\, P\in H^1(\Omega)\,,\quad\, \sigma\in H^1(\Omega) \,, \quad\, m\in H(\Curl;\Omega) \,, \quad \, \Curl m \in H^1(\Omega)\,,
\end{equation}
where the last regularity in \eqref{mainreg} follows from equation $\eqref{Main}_2$. 
Moreover, if $\bbL_c$ has a special block diagonal structure (see Corollary 
\ref{corr.reg}), then we in addition have
\begin{equation}
\label{mainreg_cor}
  m\in H^1(\Omega)\,,\quad\, \Curl P\in H^1(\Omega)\,.
\end{equation}
 The results \eqref{mainreg} are obtained by a combination of the Helmholtz decomposition for the matrix field $P$, regularity results for linear elliptic systems of elasticity-type and the classical Maxwell embedding recalled in Theorem~\ref{thm:embedding}.  The additional regularity formulated in \eqref{mainreg_cor} relies on a weighted version of the Maxwell embedding theorem, \cite{Weber81}.
\section{Background from function space theory}

\subsection{Notation, assumptions}  
\label{s:notation}

For  vectors $a,\,b\in\R^n$, we define the scalar product $\langle a,b\rangle:=\sum_{i=1}^n a_ib_i$, the Euclidean norm $\abs{a}^2:=\langle a,a\rangle$ and the dyadic product $a\otimes  b=(a_ib_j)_{i,j=1}^n\in\R^{n\times n}$, where $\R^{m\times n}$ will denote the set of real $m\times n$ matrices. 
%The identity tensor on $\R^{n\times n}$ will be denoted by $\idd$.
For  matrices $P,Q\in \R^{m\times n}$, we 
define the standard Euclidean scalar product $\langle P,Q\rangle:=\sum_{i=1}^m\sum_{j=1}^n P_{ij}Q_{ij}$ and the Frobenius-norm $\|P\|^2:=\langle P,P\rangle$. $P^T\in \R^{n\times m}$ denotes the transposition of the matrix $P\in \R^{m\times n}$ and for $P\in \R^{n\times n}$, the symmetric part of $P$
will be denoted by $\sym P=\frac{1}{2}(P+P^T)\in \Sym(n)$. 
%We recall the following scalar product rule for matrices $P\in \R^{n \times k}$, $Q\in \R^{k \times m}$ and $R\in \R^{n \times m}$: 
%\begin{align}
%\label{eq:sp_matrix}
%\langle PQ,R\rangle =\langle P, RQ^T\rangle = \langle Q, P^TR\rangle.
%\end{align}

From now on we fix $d=3$.  
Let $\Omega\subset\R^3$ be a bounded domain. 
As a minimal requirement, we assume in this paper that the boundary $\partial\Omega$ is  Lipschitz continuous, meaning that it can locally be described as the graph of a Lipschitz continuous function, see \cite{Grisvard} for a precise definition. In a similar spirit we speak of $C^1$ or $C^{1,1}$-regular boundaries. For a function $u=(u^1,\ldots, u^m)^T:\Omega\to\R^m$ (with $m\in \N$), the  differential  $\nabla u$ is given by
\begin{align*}
 \nabla u%=\Big( u^k_{,x_l}\Big)_{\substack{1\leq k\leq m,\\ 1\leq l\leq d}}
 =\begin{pmatrix}
 \nabla u^1\\
 \vdots\\
 \nabla u^m
 \end{pmatrix}
 \in \R^{m\times d}\,,
\end{align*}
with $(\nabla u^k)_{l}= \partial_{x_l} u^k$ for $1\leq k\leq m$ and $1\leq \ell\leq d$ and with $\nabla u^k\in \R^{1\times d}$. 
%,
%while for scalar functions $u:\Omega\to\R$ 
%we also define $\grad u:=(\nabla u)^T\in \R^{d\times 1}$.   
For a vector field $w:\Omega\to\R^3$, the divergence and the  curl are  given as
\begin{align*}
\dyw w=\sum_{i=1}^ 3 w^{i}_{, x_i},\qquad 
\curl w =\big(w^{3}_{,x_2}-w^{2}_{,x_3},w^{1}_{,x_3}-w^{3}_{,x_1},w^{1}_{,x_2}-w^{2}_{,x_1}\big)\,.
\end{align*}
For tensor fields $Q:\Omega\to\R^{n\times 3}$ ($n\in \N$), $\Curl Q$ and $\Dyw Q$ are defined row-wise: 
\begin{align*}
 \Curl Q=\begin{pmatrix} \curl Q^1
             \\
             \vdots\\
             \curl Q^n    
            \end{pmatrix}
\in  \R^{n\times 3}\,,\qquad 
\text{ and } 
%\begin{align*}
 \Dyw Q=\begin{pmatrix}
  \dyw Q^1
  \\ \vdots \\
   \dyw Q^n
   \end{pmatrix}\in \R^n\,,
\end{align*}
where $Q^i$ denotes the $i$-th row of $Q$. With these definitions, for $u:\Omega\to \R^m$ we have consistently $ \Curl\nabla u=0\in \R^{m \times d}$.  

%$C_0^{\infty}(\Omega)$ is the set of smooth functions with compact support in $\Omega$. 
%Additionally,  $L^2(\Omega)$ denotes the usual Lebesgue spaces of square integrable
%scalar functions, vector or tensor fields on $\Omega$ with values in $\R$, $\R^d$ or
%$\R^{m\times n}$, respectively. If needed, we will write more explicitly 
%$L^2(\Omega,\R^s)$ to indicate that we deal with $\R^s$-valued functions.
The Sobolev spaces \cite{adamssobolev,Giraultbook} used in this paper are
\begin{align*}
& H^1(\Omega)=\{u\in L^2(\Omega)\, |\,\, \nabla u\in L^2(\Omega)\}\,, \qquad \|u\|^2_{H^1(\Omega)}:=\|u\|^2_{L^2(\Omega)}+\|\nabla u\|^2_{L^2(\Omega)}\, ,\\[2ex]
&H({\curl};\Omega)=\{v\in L^2(\Omega;\R^d)\, |\,\, \curl v\in L^2(\Omega)\}\,,\qquad\|v\|^2_{H({\rm curl};\Omega)}:=\|v\|^2_{L^2(\Omega)}+\|{\rm curl}\, v\|^2_{L^2(\Omega)}\, ,\\[2ex]
&H(\dyw;\Omega)=\{v\in L^2(\Omega;\R^d)\, |\,\, \dyw v\in L^2(\Omega)\},\qquad
\|v\|^2_{H({\dyw};\Omega)}:=\|v\|^2_{L^2(\Omega)}+\|\dyw v\|^2_{L^2(\Omega)}\, ,
\end{align*}
spaces for tensor valued functions are denoted by $H(\Curl;\Omega)$ and $H(\Dyw;\Omega)$. 
Moreover, $H_0^1(\Omega)$ is the completion of $C_0^{\infty}(\Omega)$ with respect to the $H^1$-norm and $H_0({\rm curl};\Omega)$ and $H_0(\dyw;\Omega)$ are the completions of $C_0^{\infty}(\Omega)$ with respect to the $H(\curl)$-norm and the $H(\dyw)$-norm, respectively. 
By $H^{-1}(\Omega)$ we denote the dual of $H_0^1(\Omega)$. Finally we define
\begin{equation} \label{space}
\begin{split}
 H(\dyw,0;\Omega)&=\{u\in H(\dyw;\Omega)\, |\,\, \dyw u=0\}\,,\\
 H(\curl,0;\Omega)&=\{u\in H(\curl;\Omega)\, |\,\, \curl u=0\}
\end{split}
\end{equation}
and set
\begin{equation} \label{space1}
\begin{split}
H_0(\dyw,0;\Omega)&=H_0(\dyw;\Omega)\cap H(\dyw,0;\Omega)\,,\\ 
 H_0(\curl,0;\Omega)&=H_0(\curl;\Omega)\cap H(\curl,0;\Omega)\,.
\end{split}
\end{equation}
%
%  Finally, 
%\begin{equation*}
%    L^2_{{\rm loc}}(\Omega)=\{u:\Omega\rightarrow\R\, |\,\,  u\in L^2(\tilde{\Omega})\,\,%\mathrm{for\, all}\,\, \tilde{\Omega}\Subset\Omega\}\,,
%\end{equation*}
%where $\tilde{\Omega}\Subset\Omega$ means that there exists $K$ compact such that  
%$\tilde{\Omega}\subset K\subset\Omega$. This definition extends in a natural way to  the above Sobolev spaces.

\par 
{\bf Assumption A:}\hspace{1ex} We assume that the coefficient functions $\D,\Ha$ and
 $\Lc$ in \eqref{Main} are fourth order elasticity tensors from $ C^{0,1}(\overline\Omega;\Lin(\R^{3\times 3};\R^{3\times 3}))$  and are symmetric and positive definite in the following sense
\begin{enumerate}[(i)]
 \item For every $\sigma,\tau\in \Sym(3)$, $\eta_1,\eta_2\in \R^{3\times 3}$ and all $x\in \overline\Omega$:
 \begin{gather}
  \langle \mathbb{C}_e(x)\sigma,\tau\rangle = \langle \sigma,\mathbb{C}_e(x)\tau\rangle\,,
  \quad\qquad
  \langle \mathbb{C}_\text{micro}(x)\sigma,\tau\rangle = 
  \langle \sigma,\mathbb{C}_\text{micro}(x)\tau\rangle\,,
  \nonumber\\
  \langle \mathbb{L}_c(x)\eta_1,\eta_2\rangle 
  =  \langle\eta_1, \mathbb{L}_c(x)\eta_2\rangle\,.
  \label{symmetries1}
 \end{gather}
\item There exists positive constants $C_{\mathrm{e}}$, $C_{\mathrm{micro}}$ and $L_{\mathrm{c}}$ such that for all $x\in\overline\Omega$, $\sigma\in\Sym(3)$ and $\eta\in\R^{3\times 3}$:
\begin{align}
    \label{positivedef}
  \langle\D(x)\sigma,\sigma\rangle\geq C_{\mathrm{e}}|\sigma|^2 &\,,\quad\,\, \langle\Ha(x)\sigma,\sigma\rangle\geq C_{\mathrm{micro}}|\sigma|^2\,,\quad\,\,   \langle\Lc(x)\eta,\eta\rangle\geq L_{\mathrm{c}}|\eta|^2\,. 
\end{align}
\end{enumerate}
%  
%  
% \begin{enumerate}[(i)]
%     \item $\D=\D(x)=\{(c_{\mathrm{e}})_{i,j,r,s}(x)\}_{i,j,r,s=1}^3$, $\Ha=\Ha(x)=\{(c_{\mathrm{micro}})_{i,j,r,s}(x)\}_{i,j,r,s=1}^3$ and $\Lc=\Lc(x)=\{(l_{\mathrm{c}})_{i,j,r,s}(x)\}_{i,j,r,s=1}^3$ have the following symmetries
% \begin{align}
%     \label{symmetries1}
%     (c_{\mathrm{e}})_{i,j,r,s}(x)&=(c_{\mathrm{e}})_{j,i,r,s}(x)=(c_{\mathrm{e}})_{r,s,i,j}(x)\,,\\[1ex]
%      \label{symmetries2}
%      (c_{\mathrm{micro}})_{i,j,r,s}(x)&=(c_{\mathrm{micro}})_{j,i,r,s}(x)=(c_{\mathrm{micro}})_{r,s,i,j}(x)\,,\\[1ex]
%       \label{symmetries3}
%      (l_{\mathrm{c}})_{i,j,r,s}(x)&=(l_{\mathrm{c}})_{r,s,i,j}(x)\,.
% \end{align}
% \item there exist positive constants $C_{\mathrm{e}}$, $C_{\mathrm{micro}}$ and $L_{\mathrm{c}}$ such that for all $x\in\overline\Omega$, $\sigma\in\Sym(3)$ and $\eta\in\R^{3\times 3}$:
% \begin{align}
%     \label{positivedef}
%   \langle\D(x)\sigma,\sigma\rangle\geq C_{\mathrm{e}}|\sigma|^2 &\,,\quad \langle\Ha(x)\sigma,\sigma\rangle\geq C_{\mathrm{micro}}|\sigma|^2\,,\quad   \langle\Lc(x)\eta,\eta\rangle\geq L_{\mathrm{c}}|\eta|^2\,. 
% \end{align}
% \end{enumerate}

\subsection{Helmholtz decomposition, embeddings, elliptic regularity}

Based on the results from Section $3.3$ of the book \cite{Giraultbook} (Corollary $3.4$), see also  \cite[Theorem 5.3]{BaPaSc16}, the following version of the Helmholtz decomposition will be used:

\begin{theorem}[Helmholtz decomposition]
\label{thm:helmholtz-dec}
Let $\Omega\subset\R^3$ be a bounded domain with a Lipschitz boundary. Then
\[
 L^2(\Omega;\R^3)= \nabla H_0^1(\Omega)\oplus H(\dyw,0;\Omega)
\]
and hence, for every $p\in L^2(\Omega;\R^3)$ there exist unique $v\in H_0^1(\Omega)$  and 
$q\in H(\dyw,0;\Omega)$ such that 
$p=\nabla v + q$. 
\end{theorem}
%Observe that for bounded Lipschitz domains we also have
%$H(\dyw,0;\Omega)=\curl\big( H(\curl;\Omega)\big)$ (\cite[Theorem 5.3]{BaPaSc16}). 
An immediate consequence of the Helmholtz decomposition theorem is
 \begin{proposition} 
 \label{prop:helmholtz-curl}
 Let $\Omega\subset\R^3$ be a bounded domain with a Lipschitz boundary and 
let $p\in H_0(\Curl;\Omega)$ with $p=\nabla v + q$, where $v\in H_0^1(\Omega)$ and $q\in H(\Dyw,0;\Omega)$ are given according to the Helmholtz decomposition. 
 Then $\nabla v\in H_0(\Curl,0;\Omega)$ %=\Set{q\in H_0(\Curl;\Omega)}{\curl q=0}$ 
 and $q\in H(\Dyw,0;\Omega)\cap H_0(\Curl;\Omega)$.
\end{proposition}
\begin{proof}
By standard arguments it follows that $\curl\nabla v=0$ in the distributional sense, and hence $\nabla v\in H(\curl,0;\Omega)$. This implies that $\nabla v\times n$ (with $n$ the exterior unit normal vector field on $\partial\Omega$) is well defined as an element from $H^{-1/2}(\partial\Omega)$, \cite[Theorem 2.11]{Giraultbook}. 
%Girault-Raviart-book or AmroucheBernardiDaugeGirault
Moreover, for all $\phi \in C^\infty(\overline\Omega,\R^3)$ one finds
\begin{align*}
 \langle (\nabla v)\times n,\phi\rangle&=\int_\Omega\langle \nabla v, \curl\phi\rangle\dx -\int_\Omega \langle\curl(\nabla v), \phi\rangle\dx\,.
\end{align*}
The second integral vanishes since $\curl \nabla v=0$. Applying Gauss' Theorem to the first integral and taking into account that $v\big|_{\partial\Omega}=0$ implies that the first integral vanishes, as well. Hence, we finally obtain $\nabla v\times n=0$ on $\partial\Omega$ in the sense of traces, \cite[Lemma 2.4, Theorem 2.12]{Giraultbook}. Since $q=p-\nabla v$ and since $\nabla v,p\in H_0(\curl,\Omega)$ the assertion on $q$ is immediate.
\end{proof}
 The next embedding theorem is for instance proved in  \cite[Sections 3.4, 3.5]{Giraultbook}:
\begin{theorem}[Embedding Theorem]
\label{thm:embedding}
 Let $\Omega\subset \R^3$ be a bounded domain with a $C^{1,1}$-smooth boundary 
 $\partial\Omega$. Then 
 \begin{align*}
 %\label{ET1}
  H(\curl;\Omega)\cap H_0(\dyw;\Omega)\subset H^1(\Omega),\qquad\qquad 
  % \label{ET2}
  H_0(\curl;\Omega)\cap H(\dyw;\Omega)\subset H^1(\Omega)
 \end{align*}
 and there exists a constant $C>0$ such that for every $p\in  H(\curl;\Omega)\cap H_0(\dyw;\Omega)$ or\\ $p\in H_0(\curl;\Omega)\cap H(\dyw;\Omega)$ we have
 \begin{equation*}
   \norm{p}_{H^1(\Omega)}\leq C(\norm{p}_{H(\curl;\Omega)} + \norm{p}_{H(\dyw;\Omega)})\,.  
 \end{equation*}
\end{theorem}
%For the proof of Theorem \ref{thm:embedding} 
%we refer to a book \cite{Giraultbook}- Section 3.4 (for second embedding)
%and Section 3.5 (for first embedding). 
A version of this result for Lipschitz domains is for instance available in \cite{Costabel90}. 
For the previous embedding theorem there are also some versions with weights, and we cite here  Theorem 2.2 from \cite{Weber81} with $k=1$ and $\ell=0$. We assume that the weight function  $\varepsilon:\overline\Omega\to\R^{3\times 3}$ for every $x\in \overline\Omega$ is symmetric and positive definite, uniformly in $x$. 

\begin{theorem}
\label{thm:embweight}
 Let $\Omega\subset\R^3$ be a bounded domain with a $C^2$-smooth boundary and let $\varepsilon\in C^1(\overline\Omega,\R^{3\times 3})$ be a symmetric and positive definite weight function. Assume that $p:\Omega\to\R^3$ belongs to one of the following spaces:
 \begin{align}
  p\in H_0(\curl;\Omega)\quad\text{and}\quad \varepsilon\, p\in H(\dyw;\Omega)
 \end{align}
 or 
 \begin{align}
  p\in H(\curl;\Omega)\quad\text{and}\quad \varepsilon\, p\in H_0(\dyw;\Omega)\,.
 \end{align}
 Then $p\in H^1(\Omega)$ and there exists a constant $C>0$ (independent of $p$) such that 
 \begin{align*}
  \norm{p}_{H^1(\Omega)}\leq C\Big(\norm{p}_{H(\curl;(\Omega)} + \norm{\dyw (\varepsilon p)}_{L^2(\Omega)}\Big)\,.
 \end{align*}
\end{theorem}
\begin{remark}
 Assuming higher regularity on the weight function $\varepsilon$ and the smoothness of $\partial\Omega$ (i.e.\ $\varepsilon\in C^k(\overline{\Omega},\R^{3\times 3})$ and $\partial\Omega\in C^{k+1})$, Theorem 2.2 from \cite{Weber81} guarantees a corresponding higher regularity of $p$.
\end{remark}
In the proof of Theorem \ref{thm:globreg} we will decompose the microdistortion tensor $P$ as $P=\nabla q + Q$ and apply Theorem \ref{thm:embedding} to $Q$. The regularity for the displacement field $u$ and the vector $q$ then is a consequence of an elliptic regularity result that we discuss next.

Let us consider  the following auxiliary bilinear form: for $(u,\,q)\in H^1_0(\Omega;\R^{3+3})$ and $(u,\,v)\in H^1_0(\Omega;\R^{3+3})$ we define
\begin{align}
\label{auxbilinear}
\tilde{a}\Big(\begin{pmatrix} 
   u\\ q 
  \end{pmatrix},\begin{pmatrix} 
   v\\ w 
  \end{pmatrix}\Big)
  =  \int_\Omega 
  \langle \bbA(x)\begin{pmatrix} 
   \sym\nabla u\\ \sym\nabla q 
  \end{pmatrix}\,,\,
  \begin{pmatrix} 
   \sym\nabla v\\ \sym\nabla w
  \end{pmatrix}\rangle 
  \dx\,,
\end{align}
where $ \bbA:\Omega\rightarrow \big(\Lin(\Sym(3),\Sym(3))\big)^4$ is defined by the following formula
\begin{align}
\label{bbA}
 \bbA(x)=\begin{pmatrix}
       \D(x) & -\D(x)\\ - \D(x) & \D(x) + \Ha(x)
      \end{pmatrix}\,.
\end{align}
\begin{corollary}
    Let Assumption $\bf{A}_{(ii)}$ be satisfied. Then, there exists a positive constant $C_{\bbA}$ such that for all $x\in\overline\Omega$ and $\sigma=(\sigma_1,\sigma_2)\in\Sym(3)\times\Sym(3)$ we have:
\begin{equation}
    \label{positivedef1}
  \langle\bbA(x)\sigma,\sigma\rangle\geq C_{\bbA}|\sigma|^2\,. 
\end{equation}
\end{corollary}
\begin{proof}
Fix $x\in\overline\Omega$ and $\sigma=(\sigma_1,\sigma_2)\in\Sym(3)\times\Sym(3)$, then 
\begin{equation}
    \label{positivedef2}
 \langle\bbA(x)\sigma,\sigma\rangle= \langle\D(\sigma_1-\sigma_2),\sigma_1-\sigma_2\rangle+\langle\Ha(\sigma_2),\sigma_2\rangle\,.\nn
\end{equation}
Assumption $\bf{A}_{(ii)}$ implies
\begin{equation}
\begin{split}
    \label{positivedef3}
 \langle\bbA(x)\sigma,\sigma\rangle &\geq C_{\mathrm{e}}|\sigma_1-\sigma_2|^2 +C_{\mathrm{micro}}|\sigma_2|^2\geq \min\{C_{\mathrm{e}},C_{\mathrm{micro}}\}\big(|\sigma_1-\sigma_2|^2+|\sigma_2|^2\big)\nn\\[1ex]
 &\geq \frac{2}{9}\min\{C_{\mathrm{e}},C_{\mathrm{micro}}\}\big(|\sigma_1|^2+|\sigma_2|^2\big)=\frac{2}{9}\min\{C_{\mathrm{e}},C_{\mathrm{micro}}\}|\sigma|^2
\end{split}
\end{equation}
and the proof is completed.
\end{proof}
\noindent
Now, for all $(u,\,q)\in H^1_0(\Omega;\R^{3+3})$ we have 
 \begin{align}
\label{bicoerciv}
\tilde{a}\Big(\begin{pmatrix} 
   u\\ q 
  \end{pmatrix},\begin{pmatrix} 
   u\\ q 
  \end{pmatrix}\Big)&\geq C_{\bbA}\big(\|\sym\nabla u\|^2_{L^2(\Omega)}+\|\sym\nabla q\|^2_{L^2(\Omega)}\big)\nn\\[1ex]
  &\geq  C_{\bbA}\,C_{K}\big(\|u\|^2_{H^1_0(\Omega)}+\|q\|^2_{H^1_0(\Omega)}\big)\,,
\end{align}
where the constant $C_{K}$ is a constant resulting from the standard Korn's inequality \cite{neff_2006}. This shows that    the bilinear form \eqref{auxbilinear} is coercive on the space  $H^1_0(\Omega;\R^{3+3})$. The form \eqref{auxbilinear} defines the following auxiliary problem: for $(F_1,F_2)\in L^2(\Omega;\R^{3+3})$ find  $(u,\,q)\in H^1_0(\Omega;\R^{3+3})$ with
\begin{align}
\label{auxproblem}
\tilde{a}\big(\left(\begin{smallmatrix} 
   u\\ q 
  \end{smallmatrix}\right),\left(\begin{smallmatrix} 
   v\\ w 
  \end{smallmatrix}\right)\big)=  \int_\Omega 
  \langle
  \left(\begin{smallmatrix} 
   F_1\\ F_2 
  \end{smallmatrix}\right),\left(\begin{smallmatrix} 
   v\\ w 
  \end{smallmatrix}\right) \rangle
  \dx 
\end{align}
for all $(v,\,w)\in H^1_0(\Omega;\R^{3+3})$. Modifying the results concerning the regularity of  elliptic partial differential equations, it would be possible to obtain the existence of a solution for system \eqref{auxproblem} with the regularity $(u,\,q)\in H^1_0(\Omega;\R^{3+3})\cap H^2(\Omega;\R^{3+3})$ (see for example \cite[Theorem 9.15,  Section 9.6]{GiTr86}). However, system \eqref{auxproblem} fits perfectly into the class  considered in \cite{NK08}.  There, the global regularity of weak solutions to a quasilinear elliptic system with a rank-one-monotone nonlinearity was investigated. 
As an application of the result from \cite{NK08} we obtain  
\begin{lemma}
\label{thm:linellreg}
Let $\Omega\subset \R^3$ be a bounded domain with a $C^{1,1}$-smooth boundary $\partial\Omega$. Let furthermore  $(F_1,F_2)\in L^2(\Omega;\R^{3+3})$ and  $\D,\Ha\in C^{0,1}(\overline\Omega; \Lin(\R^{3\times 3};\R^{3\times 3}))$. Then the problem \eqref{auxproblem} has a unique solution $(u,\,q)\in H^1_0(\Omega;\R^{3+3})\cap H^2(\Omega;\R^{3+3})$.
 \end{lemma}
\begin{proof}
Coercivity \eqref{bicoerciv} of the bilinear form \eqref{bilinear} and the Lax-Milgram Theorem imply the existence of exactly one weak solution  $(u,\,q)\in H^1_0(\Omega;\R^{3+3})$. In order to prove higher regularity of this solution we will use  Theorem $5.2$ of \cite{NK08}.

Let us introduce  $\mathbb{B}:\overline\Omega\times\R^{(3+3)\times 3}\rightarrow \R^{(3+3)\times 3}$ as the unique $x$-dependent linear mapping satisfying 
\begin{align}
\label{operatorB}
 \langle \mathbb{B}\Big(x,
 \left(\begin{smallmatrix} 
   A_1\\ A_2 
  \end{smallmatrix}\right)
  \Big),
  \left(\begin{smallmatrix} 
   B_1\\ B_2 
  \end{smallmatrix}\right)\rangle 
  =\langle \bbA(x)
  \left(\begin{smallmatrix} 
   \sym A_1\\ \sym A_2  
  \end{smallmatrix}\right)\,,\,
  \left(\begin{smallmatrix} 
   \sym B_1\\ \sym B_2
  \end{smallmatrix}\right)
\rangle 
\end{align}
for all $A_1,A_2,B_1,B_2\in \R^{3\times 3}$ and $x\in \overline\Omega$. 
Then for every $x\in\overline{\Omega}$, $A=(A_1,A_2)\in\R^{(3+3)\times 3}$, $\xi=(\xi_1,\xi_2)\in\R^{3+3}$ and $\eta\in\R^3$ we have thanks to the positive definiteness of  $\bbA$ 
\begin{equation}
\begin{split}
\label{rankone}
\big\langle\,\mathbb{B}\Big(x,\left(\begin{smallmatrix} 
   A_1\\ A_2 
  \end{smallmatrix}\right)
  +\xi\otimes\eta\Big)&-\mathbb{B}\Big(x,
  \left(\begin{smallmatrix} 
   A_1\\ A_2 
  \end{smallmatrix}
  \right)
  \Big),\xi\otimes\eta\,\big\rangle
  %&=\big\langle\,\mathbb{B}\big(x,\xi\otimes\eta\big),\xi\otimes\eta\,\big\rangle\nn \\
  %=\big\langle\,\bbA(x)
  %\begin{pmatrix} 
  % \sym(\xi_1\otimes\eta)
  % \\
  % \sym(\xi_2\otimes\eta)
%\end{pmatrix}
%   ,
%    \begin{pmatrix} 
%  \sym( \xi_1\otimes\eta)
%   \\
%   \sym(\xi_2\otimes\eta)
%\end{pmatrix}
%   \,\big\rangle
%\\   %\nonumber
   \geq\\[1ex] 
   &\geq C_{\bbA}
   \big(\norm{\sym(\xi_1\otimes \eta)}^2 + \norm{\sym(\xi_2\otimes \eta)}^2\big)\,.
   \end{split}
\end{equation}
 Since $\norm{\sym(\xi_i\otimes \eta)}^2=\frac12(\norm{\xi_i}^2\norm{\eta}^2 + \norm{\langle \xi_i,\eta\rangle}^2)$, this ultimately implies that $\mathbb{B}$ is strongly rank-one monotone/satisfies the Legendre-Hadamard condition. 
 Due to the Lipschitz continuity of $\D$ and $\Ha$ the remaining assumptions of 
 Theorem $5.2$ of \cite{NK08} can easily be verified.
% there exist  constants $L_1$, $L_2>0$ such
%that for every $x,\, x_i\in\Omega$, $A,\, A^i\in\R^{(3+3)\times 3}$ 
%\begin{align}
%\label{continuity}
%\|\mathbb{B}(x_1,A)-\mathbb{B}(x_2,A)\|&\leq L_1\|x_1-x_2\|\|A\|\,,\nn\\
%\|\mathbb{B}(x,A^1)-\mathbb{B}(x,A^2)\|&\leq L_2\|A^1-A^2\|\,,\nn\\
%\mathbb{B}(x,0)&=0\,.
%\end{align}
%%With properties \eqref{bicoerciv}, \eqref{rankone} and \eqref{continuity}
Hence,   \cite[Theorem 5.2]{NK08} implies  $(u,\,q)\in H^2(\Omega;\R^{3+3})$.
\end{proof}

\section{Weak formulation of the relaxed micromorphic model}
For $u,\,v\in H^1_0(\Omega,\R^3)$ and $P,\,W\in H_0(\Curl;\Omega)$ the following bilinear form is associated with the system \eqref{Main}
\begin{align}
\label{bilinear}
a\big((u,P),(v,W)\big)&=\int_{\Omega}\Big(\big\langle \D\sym(\nabla u-P),\sym(\nabla v-W)\big\rangle +\big\langle\Ha\sym P,\sym W\big\rangle\\
&\qquad +
\big\langle\Lc\Curl P,\Curl W\big\rangle\Big)\,\di x
\nonumber 
\\[1ex]
&\equiv \int_\Omega \langle \bbA 
\begin{pmatrix} 
\sym\nabla u\\
\sym P
\end{pmatrix}
\,,\,
\begin{pmatrix} 
\sym \nabla v\\
\sym W
\end{pmatrix}\rangle
+ \big\langle\Lc\Curl P,\Curl W\big\rangle\Big)\,\di x
\nonumber 
\end{align}
where the tensor $\bbA$ is defined in \eqref{bbA}. Here, homogeneous boundary conditions $u\big|_{\partial\Omega}=0$ and $(P\times n)\big|_{\partial\Omega}=0$ are considered. 
\begin{theorem}[Existence of weak solutions]
\label{thm:weakex}
Let $\Omega\subset\R^3$ be a bounded domain with a Lipschitz boundary and 
assume that $\D,\Ha,\Lc\in L^\infty(\Omega;\Lin(\R^{3\times 3};\R^{3\times 3}))$ comply  with the symmetry and positivity properties formulated in \eqref{symmetries1} - \eqref{positivedef}. Then for every $f\in H^{-1}(\Omega)$ and $M\in \big(H_0(\Curl;\Omega)\big)^{\ast}$ there exists a unique pair $(u,P)\in H_0^1(\Omega)\times H_0(\Curl;\Omega)$ such that 
\begin{align}
\label{weakform}
\forall (v,W)\in H_0^1(\Omega)\times H_0(\Curl;\Omega):\qquad \,\,
a\big((u,P),(v,W)\big)=\int_{\Omega}\langle f,v\rangle + \langle M,W\rangle\,\di x\,.
\end{align}
\end{theorem}
\begin{proof}
For a bounded domain $\Omega\subset\R^3$ with Lipschitz boundary $\partial\Omega$ the incompatible Korn's inequality  \cite{OptimalKMSIneq,lewintan_neff_2021,NeffLP,neff2012poincare} implies that there is a constant $\tilde{c}>0$ such that 
\begin{equation}
\label{coercive1}
\|P\|^2_{L^2(\Omega)}\leq \tilde{c}\,\big( \|\sym  P\|^2_{L^2(\Omega)}+\|\Curl P\|^2_{L^2(\Omega)}\big)
\end{equation}
for all $P\in {\rm H}_0(\Curl;\Omega)$. Positive definiteness of the tensors $\bbA$ and $\Lc$ entail
\begin{equation}
\label{Coercive}
a\big((u,P),(u,P)\big)\geq C_{\bbA}(\|\sym \nabla u\|^2_{L^2(\Omega)}+\|\sym P\|^2_{L^2(\Omega)})+L_{\mathrm{c}}\|\Curl P\|^2_{L^2(\Omega)}\,.
\end{equation}
Hence, by \eqref{coercive1} and Korn's inequality the bilinear form \eqref{bilinear} is coercive on $H_0^1(\Omega)\times H_0(\Curl;\Omega)$ and  the Lax-Milgram Theorem finishes the proof.
%in the sense that there exists positive constant $C>0$ such that for all $u\in {\rm H}%^1_0(\Omega)$ and $P\in {\rm H}_0(\Curl;\Omega)$ it holds
%\begin{equation}
%\label{Coercive1}
%a\big((u,P),(u,P)\big)\geq C(\| u\|^2_{H^1(\Omega)}+\|P\|^2_{L^2(\Omega)}+\|\Curl P\|%^2_{L^2(\Omega)})\,.
%\end{equation}
\end{proof}
Thanks to the Helmholtz decomposition, weak solutions can equivalently be characterized as follows:
\begin{lemma}
Let the assumptions of Theorem \ref{thm:weakex} be satisfied, $f\in H^{-1}(\Omega)$ and $M\in \big(H_0(\Curl;\Omega)\big)^{\ast}$. Let furthermore $(u,P)\in H_0^1(\Omega)\times H_0(\Curl;\Omega)$ and let $(q,Q)\in H_0^1(\Omega;\R^3)\times H(\Dyw,0;\Omega)$ such that $P=\nabla q + Q$. Then the following (a) and (b) are equivalent:
\begin{itemize}
 \item[(a)] $(u,P)$ is a weak solution of \eqref{Main} in the sense of \eqref{weakform}.
 \item[(b)] For all $(v,W)\in H_0^1(\Omega)\times H_0(\Curl;\Omega)$ the triple $(u,q,Q)$ satisfies 
 \begin{align}
 \label{weakformhelmholtz}
  \int_\Omega 
  \langle \bbA & 
  \begin{pmatrix} 
   \sym\nabla u\\ \sym\nabla q + \sym Q
  \end{pmatrix}\,,\,
  \begin{pmatrix} 
   \sym\nabla v\\ \sym W
  \end{pmatrix}
  \rangle
  \dx\nn\\[1ex]&
  + \int_\Omega \langle\Lc \Curl Q,\Curl W\rangle\dx= \int_{\Omega}\langle f,v\rangle + \langle M,W\rangle\,\di x\,. 
 \end{align}
\end{itemize}
\end{lemma}
\begin{proof}
%This is an immediate consequence of the identity $\Curl\nabla q=0$ for $q\in H^1(\Omega)$.
$(a)\implies (b):$\hspace{1ex} Assume that $(u,P)\in H_0^1(\Omega)\times H_0(\Curl;\Omega)$ is a weak solution of \eqref{Main} in the sense of \eqref{weakform}. Then Theorem \ref{thm:helmholtz-dec} implies that for $i=1,2,3$ there exists unique $q_i\in H_0^1(\Omega)$  and $Q_i\in H_0(\Dyw,0;\Omega)$ such that $P_i=\nabla q_i + Q_i$, where $P_i$ denotes the rows of the matrix $P$. Inserting $P=(\nabla q_1 + Q_1, \nabla q_2 + Q_2, \nabla q_3 + Q_3)^T$ into \eqref{weakform} we obtain \eqref{weakformhelmholtz}, where $Q=(Q_1, Q_2, Q_3)^T$\,.\\
$(b)\implies (a):$\hspace{1ex} Let for all $(v,W)\in H_0^1(\Omega)\times H_0(\Curl;\Omega)$ the triple $(u,q,Q)\in H_0^1(\Omega;\R^3)\times H_0^1(\Omega;\R^3)\times H(\Dyw,0;\Omega)$ satisfy \eqref{weakformhelmholtz}. Then \eqref{weakformhelmholtz} can be written in the form
\begin{align}
 \label{weakformhelmholtz1}
  \int_\Omega \langle\bbA &
  \begin{pmatrix} 
   \sym\nabla u\\ \sym\nabla q + \sym Q
  \end{pmatrix},
  \begin{pmatrix} 
   \sym\nabla v\\ \sym W
  \end{pmatrix}\rangle
  \dx\nn\\[1ex]
  &+ \int_\Omega \langle\Lc \Curl (\nabla q+Q),\Curl W\rangle\,\dx= \int_{\Omega}\langle f,v\rangle + \langle M,W\rangle\,\di x
 \end{align}
and the function $(u,\nabla q+Q)$ satisfies \eqref{weakform} for all $(v,W)\in H_0^1(\Omega)\times H_0(\Curl;\Omega)$. Uniqueness of a weak solution of the problem \eqref{Main} implies that $P=\nabla q+Q$.
\end{proof}
\section{Global regularity on smooth domains}
\label{sec:globreg}
The aim of this section is to prove the following regularity theorem
\begin{theorem}
\label{thm:globreg}
 Let $\Omega\subset \R^3$ be a bounded domain with a $C^{1,1}$-smooth boundary. Moreover, in addition to the assumptions of Theorem \ref{thm:weakex} let $\D,\Ha,\Lc\in C^{0,1}(\overline\Omega; \Lin(\R^{3\times 3};\R^{3\times 3}))$. 
 Finally, we assume that $f\in L^2(\Omega)$ and $M\in H(\Dyw;\Omega)$. Then for every weak solution   $(u,P)\in H_0^1(\Omega)\times H_0(\Curl;\Omega)$ we have 
 \begin{align}
  u\in H^2(\Omega)\,,\qquad\qquad P\in H^1(\Omega)\,,\qquad\qquad \Lc\Curl P\in  H(\Curl;\Omega)
 \end{align}
 and there exists a constant $C>0$ (independent of $f$ and $M$) such that 
 \begin{align}
  \norm{u}_{H^2(\Omega)} + \norm{P}_{H^1(\Omega)} + \norm{\mathbb{L}_c\Curl P}_{H(\Curl;\Omega)} \leq C(\norm{f}_{L^2(\Omega)} + \norm{M}_{H(\Dyw;\Omega)})\,.
 \end{align}
\end{theorem}
The proof relies on the Helmholtz decomposition of $P$, the embedding Theorem \ref{thm:embedding} and Theorem \ref{thm:linellreg} about the global regularity for the  auxiliary problem \eqref{auxproblem}.
\begin{proof}
 Let $(u,P)\in H_0^1(\Omega)\times H_0(\Curl;\Omega)$ satisfy \eqref{weakform} with $f\in L^2(\Omega)$ and $M\in H(\Dyw;\Omega)$.   
 
 We first show that $(u,P)\in H^2(\Omega)\times H^1(\Omega)$. 
 Let  $P=\nabla q + Q$, where $q\in H_0^1(\Omega;\R^3)$ and $Q\in H(\Dyw,0;\Omega)$ are given according to the Helmholtz decomposition. By Proposition \ref{prop:helmholtz-curl} it follows that $Q\in  H(\Dyw,0;\Omega)\cap H_0(\Curl;\Omega)$ and thanks to the assumed regularity of $\partial \Omega$, Theorem \ref{thm:embedding} implies that $Q\in H^1(\Omega)$. Next,  choosing $W=\nabla w$ for $w\in C_0^\infty(\Omega;\R^d)$,  the weak form \eqref{weakformhelmholtz} in combination with a density argument implies that for all  $v,w\in H_0^1(\Omega)$ we have
 \begin{equation}
 \label{4.3}
 \begin{split}
  \int_\Omega \langle\bbA & \begin{pmatrix} 
   \sym\nabla u\\ \sym\nabla q 
  \end{pmatrix}\,,\,
  \begin{pmatrix} 
   \sym\nabla v\\ \sym\nabla w
  \end{pmatrix}
  \rangle
  \dx\\[1ex]=&\int_\Omega \langle\D \sym Q,\sym\nabla v\rangle\dx 
  -\int_\Omega \langle(\D +\Ha)\sym Q,\sym\nabla w\rangle\dx\\[1ex]
  &+\int_\Omega \langle f, v\rangle + \langle M,\nabla w\rangle \dx\,. 
 \end{split}
  \end{equation}
Since $Q\in H^1(\Omega)$ and $M\in H(\Dyw;\Omega)$, by partial integration the right hand side of \eqref{4.3}  can be rewritten as $\int_\Omega \langle F_1, v\rangle + \langle F_2, w\rangle\dx$ with  functions $F_1,F_2\in L^2(\Omega)$.  Theorem \ref{thm:linellreg} implies that $u,q\in H^2(\Omega)$ and hence $P=\nabla q +Q\in H^1(\Omega)$. 

Let us next 
 choose $v=0$ and $W\in C_0^\infty(\Omega)$ in \eqref{bilinear}. Rearranging the terms  we find that 
\begin{align*}
 \int_\Omega \langle \Lc\Curl P,\Curl W\rangle\dx = \int_\Omega \langle M,W\rangle\dx 
 + 
 \int_\Omega \langle\D \sym\nabla u - (\D + \Ha )\sym P\,,\,\sym W\rangle\dx 
\end{align*}
which implies that $\Curl(\Lc\Curl P)\in L^2(\Omega)$ and $\Lc\Curl P\in H(\Curl;\Omega)$. 
\end{proof}

If we additionally assume  that $\mathbb{L}_c$ has a block-diagonal structure, we may also achieve $\Curl P\in H^1(\Omega)$ by applying the weighted embedding Theorem \ref{thm:embweight}.

\begin{corollary}
\label{corr.reg}
 In addition to the assumptions of Theorem \ref{thm:globreg} let $\mathbb{L}_c\in C^1(\overline\Omega;\Lin(\R^{3\times 3},\R^{3\times 3}))$ be of block diagonal structure, meaning that there exist $\mathbb{L}_i\in C^1(\overline\Omega;\R^{3\times 3})$, $1\leq i\leq 3$, such that for every $W\in \R^{3\times 3}$ we have 
 $  (\mathbb{L}_c W)_{i\text{-th row}}= \mathbb{L}_i (W_{i\text{-th row}})$. 
 Then $\mathbb{L}_c\Curl P\in H^1(\Omega)$ and  $\Curl P\in H^1(\Omega)$.
 \end{corollary}

\begin{proof}
 We focus on the $i$-th row $P^i$  of $P$. Let $\varepsilon=(\mathbb{L}_i)^{-1}$. Then 
 $\varepsilon\in C^1(\overline\Omega;\R^{3\times 3})$ and $\varepsilon(x)$ is symmetric and uniformly positive definite with respect to $x\in \overline \Omega$. 
 
 Clearly, $\varepsilon \mathbb{L}_i \curl P^i=\curl P^i \in H(\dyw,0;\Omega)$. Moreover, for the normal trace we find: for every $\phi\in C^\infty(\overline\Omega;\R)$
 \begin{align}
  \langle (\varepsilon \mathbb{L}_i \curl P^i)\cdot n,\phi\rangle_{\partial\Omega} &= 
  \langle (\curl P^i) \cdot n,\phi\rangle_{\partial\Omega}
  =\int_\Omega \langle\curl P^i, \nabla \phi\rangle\dx + \int_\Omega \langle\phi, \dyw(\curl P^i)\rangle\,\dx\nn\\
 & =\langle  P^i\times n,\phi  \rangle_{\partial\Omega}
 + \int_\Omega \langle P^i,\curl\nabla \phi\rangle\dx + \int_\Omega \langle\phi, \dyw(\curl P^i)\rangle\,\dx, 
 \end{align}
where we applied the corresponding Green's formulae \cite[Theorem 2.11, Theorem 2.5]{Giraultbook}. Since all the  terms in the last line are zero, we  obtain  
\[
 \mathbb{L}_i \curl P^i\in H(\curl;\Omega)\quad\text{and}\quad 
 \varepsilon \mathbb{L}_i \curl P^i \in H_0(\dyw,0;\Omega)\,.
\]
The weigthed embedding Theorem \ref{thm:embweight}  implies $\mathbb{L}_i \curl P^i\in H^1(\Omega)$, and since $\varepsilon=\mathbb{L}_i^{-1}$ is a multiplyer on $H^1(\Omega)$, we finally obtain $\curl P^i\in H^1(\Omega)$.% 
%consider the product rule in order to obtain the multiplyer-property
\end{proof}
\begin{remark}
The previous result may be applied to the simple uni-constant isotropic curvature case $L^2_{\mathrm{c}}\|\Curl P\|^2$.
\end{remark}
\begin{remark}
It is clear that the same higher regularity result can be established for the linear Cosserat model \cite{Ghiba2022CosseratME}.
\end{remark}
\noindent
{\it Open Problem:}\hspace{1ex} It would be interesting to establish higher regularity also for $m=\mathbb{L}_c\Curl P:=\widehat{\mathbb{L}}_c\sym\Curl P$ with $\widehat{\mathbb{L}}_c$ positive definite on symmetric arguments. However, the simple extension of the present argument fails since our argument relies on the a-priori information that $P$ belongs to $H(\Curl;\Omega)$. For the more general model involving $\widehat{\mathbb{L}}_c\sym\Curl P$ we have weak solutions in $H(\sym\Curl;\Omega)$, only, meaning that it is not clear whether in this case, $\Curl P\in L^2(\Omega)$. 

\section{Global regularity for a gauge-invariant incompatible elasticity model}
The method presented above for obtaining regularity of solution for the static relaxed micromorphic model can be directly applied to the following gauge-invariant incompatible elasticity model \cite{Lazarqauge,Lazar_Neff_gauge}
\begin{align}
\label{ecst}
0=-\Curl[\Lc\Curl e] - \D\sym e-\Dc\skyw e+M\,,
\end{align}
where the unknown function is the non-symmetric incompatible elastic distortion $e:\Omega\rightarrow \R^{3\times3}$ while $M:\Omega\rightarrow \R^{3\times 3}$ is a given body moment tensor. The constitutive tensors $\D$, $\Lc$ are positive definite fourth order tensors (fulfilling the Assumption A) while  $\Dc:\so(3)\rightarrow\so(3)$ is positive semi-definite. The system \eqref{ecst} is considered with homogeneous tangential boundary conditions, i.e.
\begin{align} \label{bc0}
{e}_i({x})\times\,n(x) =0 \ \ \
\ \textrm{for} \ \ \ {x}\in\partial \Omega\,,
\end{align}
where $\times$ denotes the vector product, $n$ is the unit outward normal vector at the surface $\partial\Omega$, $e_i$ ($i=1,2,3$) are the rows of the tensor $e$. Problem \eqref{ecst} generalises the time-harmonic Maxwell-type eigenvalue problem \cite{Alberi_maxwel,YIN_gauge} from the vectorial to the tensorial setting. On the other hand, equation \eqref{ecst} corresponds to the second equation of \eqref{Main} upon setting $\Ha\equiv 0$, assuming $\Dc\equiv 0$, identifying the elastic distortion $e$ with $e=\nabla u-P$ and observing that 
\begin{equation}
    \label{invarianc0}
    -\Curl\Lc\Curl P=\Curl\Lc\Curl(-P)=\Curl\Lc\Curl(\nabla u- P)=\Curl\Lc\Curl e\,.
\end{equation}
Smooth solutions of \eqref{ecst} satisfy the balance of linear momentum equation
\begin{align}
\label{ecst1}
\Dyw (\D\sym e+\Dc\skyw e)=\Dyw M\,.
\end{align}
Gauge-invariance means here that the solution $e$ is invariant under
\begin{equation}
\label{invariance}
    \nabla u \rightarrow \nabla u+\nabla\tau\,,\qquad\qquad P\rightarrow P+\nabla\tau
\end{equation}
which invariance is only possible since $\Ha\equiv 0$ ($\Ha> 0$ breaks the gauge-invariance). Here $\tau$ is a space-dependent (or local) translation vector. The elastic energy can then be expressed as
\begin{equation}
\label{Variational1}
\int_{\Omega}\big\langle \D\sym e,\sym e\big\rangle
+\big\langle\Lc\Curl e,\Curl e\big\rangle-\langle M,e\rangle\,\di x
\end{equation}
in which the first term accounts for the energy due to elastic distortion, the second term takes into account the energy due to incompatibility in the presence of dislocations and the last term is representing the forcing.

For $e,\,v\in H_0(\Curl;\Omega)$ the following bilinear form is associated with the system \eqref{ecst}
\begin{equation}
\label{bilinear1}
\begin{split}
b(e,v)=\int_{\Omega}\big\langle\Lc\Curl e,\Curl v\big\rangle+\big\langle\D\sym e,\sym v\big\rangle+\big\langle\Dc\skyw e,\skyw v\big\rangle
\,\di x\,.
\end{split}
\end{equation}
Let us assume that $M\in H(\Dyw;\Omega)$. Coercivity of the bilinear form \eqref{bilinear1} (the generalized incompatible Korn's inequality \eqref{coercive1}) and the Lax-Milgram Theorem imply the existence of exactly one weak solution $e\in H_0(\Curl;\Omega)$ of the system \eqref{ecst}. The Helmholtz decomposition yields that  $e=\nabla q + Q$, where $q\in H_0^1(\Omega;\R^3)$ and $Q\in H(\Dyw,0;\Omega)$. By Proposition \ref{prop:helmholtz-curl} it follows that $Q\in  H(\Dyw,0;\Omega)\cap H_0(\Curl;\Omega)$ and Theorem \ref{thm:embedding} implies that $Q\in H^1(\Omega)$. Inserting the decomposed form of the tensor $e$ into the weak form of system \eqref{ecst}, we obtain 
\begin{align}
 \label{weakformhelmholtz2}
  \int_\Omega \langle& \D\sym(\nabla q + Q), \sym W\rangle
  \dx+\int_\Omega \langle \Dc\skyw(\nabla q + Q), \skyw W\rangle
  \dx\nn\\[1ex]
  &+ \int_\Omega \langle\Lc \Curl (\nabla q+Q),\Curl W\rangle\,\dx= \int_{\Omega} \langle M,W\rangle\,\di x
 \end{align}
for all $W\in H_0(\Curl;\Omega)$. Again, choosing $W=\nabla w$ for $w\in C_0^\infty(\Omega;\R^3)$, the weak form \eqref{weakformhelmholtz2} in combination with a density argument implies that for all  $w\in H_0^1(\Omega)$ we have
\begin{align}
 \label{weakformhelmholtz3}
  \int_\Omega \langle& \D\sym\nabla q, \sym \nabla w\rangle
  \dx+\int_\Omega \langle \Dc\skyw\nabla q, \skyw \nabla w\rangle
  \dx\nn\\[1ex]
  &= \int_{\Omega} \langle M,\nabla w\rangle\,\di x-\int_\Omega \langle \D\sym Q, \sym \nabla w\rangle
  \dx-\int_\Omega \langle \Dc\skyw Q, \skyw \nabla w\rangle\,\dx
 \end{align}
Since $Q\in H^1(\Omega)$ and $M\in H(\Dyw;\Omega)$, by partial integration the right hand side of \eqref{4.3}  can be rewritten as $\int_\Omega \langle F_2, w\rangle\dx$ with  functions $F_2\in L^2(\Omega)$. Note that the auxiliary problem \eqref{weakformhelmholtz3} obtained this time is a problem from standard linear elasticity. Thus, in this case we do not need to go through Lemma \ref{thm:linellreg}: just from the standard theory of regularity in the linear elasticity we obtain that $q\in H^2(\Omega)$ (\cite{valent}). The result is that on smooth domains and with smooth coefficients the weak solution of \eqref{ecst} is more regular and satisfies 
\begin{equation}
\label{mainreg1}
e\in H^1(\Omega) \,, \quad\, \Lc\Curl e\in H(\Curl;\Omega) \,, \quad \, \Curl(\Lc\Curl e) \in H^1(\Omega)\,,
\end{equation}
where the last regularity in \eqref{mainreg1} follows from equation $\eqref{ecst}$. Moreover, if $\Lc$ has a special
block diagonal structure, then we have
\begin{equation}
\label{mainreg2}
\Lc\Curl e\in H^1(\Omega) \,, \quad \, \Curl e \in H^1(\Omega)\,.
\end{equation}

\subsubsection*{Acknowledgment} \footnotesize{
The authors wish to thank Dirk Pauly (TU Dresden) for inspiring discussions on this subject dating to 2013. Patrizio Neff and Dorothee Knees acknowledge  support  in the framework of the Priority Programme SPP 2256 "Variational Methods for Predicting Complex Phenomena in Engineering Structures and Materials" funded by the Deutsche Forschungsgemeinschaft (DFG, German research foundation): P.\ Neff within the project "A variational scale-dependent transition scheme - from Cauchy elasticity to the relaxed micromorphic continuum" (Project-ID 440935806), D.\ Knees within the project  "Rate-independent systems in solid mechanics: physical properties, mathematical analysis, efficient numerical algorithms"  (Project-ID 441222077). D.K.\ and P.N.\  enjoyed the welcoming athmosphere at the Hausdorff Research Institute for Mathematics, Bonn, funded by the Deutsche Forschungsgemeinschaft (DFG, German Research Foundation) under Germany's Excellence Strategy -- EXC-2047/1 -- 390685813.
}

\bibliographystyle{plain}
\begin{footnotesize}

\end{footnotesize}

\end{document}